\documentclass[12pt]{article}

\usepackage{amsmath,amssymb,amstext}
\usepackage{amsfonts}

\usepackage[latin1]{inputenc}

\usepackage{setspace}
\def\nbOne{{\mathchoice {\rm 1\mskip-4mu l} {\rm 1\mskip-4mu l} {\rm
1\mskip-4.5mu l} {\rm 1\mskip-5mu l}}}
\def\S{\mathcal{S}^{\downarrow}}

\newenvironment{disarray}{\everymath{\displaystyle\everymath{}}\array}
{\endarray}
\newtheorem{theo}{Theorem}

\newtheorem{cor}[theo]{Corollary}
\newtheorem{lm}{Lemma}

\newtheorem{rk}{Remark}

\newtheorem{df}{Definition}
\newenvironment{dem}{\textbf{Proof.}}{\flushright$\blacksquare$\\}

\newcommand{\EE}{\mathbb{E}}

\newcommand{\NN}{\mathbb{N}}

\newcommand{\x}{\mathbf{x}}
\newcommand{\X}{\mathbf{X}}
\newcommand{\s}{\mathbf{s}}
\newcommand{\K}{\mathbf{K}}
\newcommand{\E}{{\cal E}}
\newcommand{\ba}{{\widehat \alpha}}
\newcommand{\bb}{{\widehat \beta}}
\newcommand{\cc}{{\widehat C}}
\newcommand{\cd}{{\widehat D}}
\newcommand{\ca}{{\widehat A}}
\newcommand{\hp}{{\widehat \phi}}
\newcommand{\ttp}{{\widetilde \phi}}
\newcommand{\hpr}{{\widehat \rho}}
\newcommand{\hL}{{\widehat L}}
\newcommand{\hX}{{\widehat{\X}}}
\newcommand{\hx}{{\widehat{\x}}}
\newcommand{\hvi}{{\widehat{\varphi}}}
\newcommand{\hnu}{{\widehat{\nu}}}
\newcommand{\hmu}{{\widehat{\mu}}}
\newcommand{\hPs}{{\widehat{\Psi}}}
\newcommand{\hps}{{\widehat{\psi}}}
\newcommand{\mE}{{\mathsf E}}
\newcommand{\mmE}{\widehat {\mE}}
\begin{document}

\title{Energy efficiency of  consecutive fragmentation processes}

\author{Joaqu\'{\i}n Fontbona, Nathalie Krell, Servet Mart\'{\i}nez}


\maketitle

\begin{abstract}

We present a first study on the  energy required to reduce a unit
mass fragment by  consecutively using several devices, as it
happens in the mining industry. Two devices are considered, which
we represent as different stochastic fragmentation processes.
Following the self-similar energy model introduced by Bertoin and
Mart\'\i nez \cite{bermar}, we compute the average energy required
to attain a size $\eta_0$ with this two-device procedure. We then
asymptotically compare, as $\eta_0$ goes to $0$ or $1$, its energy
requirement with that  of individual fragmentation processes. In
particular, we show that for certain range of parameters of the
fragmentation processes and of their energy cost-functions, the
consecutive use of two devices can be asymptotically  more
efficient than using each of them separately, or conversely.

\end{abstract}

\medskip

Keywords: fragmentation process, fragmentation energy,
subordinators, Laplace exponents.

\medskip

Mathematics Subject Classification (2000): Primary 60J85,
Secondary 60J80.

\section {Introduction}

The present work is motivated by the mining industry, where
mechanical devices are used to break rocks in order to liberate
the metal contained in them. This fragmentation procedure is
carried out in a series of steps (the first of them being
blasting, followed then by  crushers, grinders or mills) until
fragments attain a sufficiently small size for the mining
purposes. One of the problems that faces the mining industry is to
minimize the total amount of energy consumed in this process. To
be more precise, at each intermediate step, material is broken by
a repetitive mechanism until  particles can go across a
classifying-grid leading to the next step. The output sizes are
known to be not optimal in terms of the global energy cost.
Moreover, since crushers or mills are large and hardly replaceable
machines,  those output sizes are in practice one of the few
parameters on which a decision can be made.

In an idealized setting,  the problem might be posed as follows:
suppose that  a unit-size fragment is to be reduced into fragments
of sizes smaller than a fixed threshold $\eta_0\in (0,1]$, by
passing consecutively through two different
 fragmentation mechanisms (for instance the first one
could be constituted by the crushers and the second one by the
mills). In this ``two-step'' fragmentation procedure,  each mass
fragment evolves in the first fragmentation mechanism until  it
first becomes smaller than $\eta\in (\eta_0,1]$, at which moment
it immediately enters the second mechanism. Then, the fragment
continues to evolve until the first instant it becomes smaller
than $\eta_0$, when it finally  exits the system. The central
question is:

\smallskip

(*) \ what is the optimal choice for the intermediate threshold
$\eta$?

\smallskip

To formulate this problem we shall model each fragmentation
mechanism by a continuous-time random fragmentation process, in
which particles break independently of each other (branching
property) and in a self-similar way. (For recent a account and
developments on the mathematical theory of fragmentation
processes, we refer to Bertoin \cite{ber}.) The self-similarity
hypothesis agrees with observations made by the mining industry;
see e.g. \cite{Charles}. In particular, it is reasonable to assume
that the energy required to break a block of size $s$ into a set
of smaller blocks of sizes $(s_1,s_2,...)$ is of the form $s^\beta
\varphi(s_1/s,s_2/s,\ldots)$, where $\varphi$ is a cost function
and $\beta >0$ a fixed parameter. For example, in the so-called
potential case, one has
$\varphi(s_1,s_2,\ldots)=\sum_{n=1}^{\infty} s_n^\beta-1$, which
corresponds to the law of Charles, Walker and Bond \cite{Charles}.

Within that mathematical framework, the asymptotic behavior of the
energy required by a single fragmentation process in order that
all fragments  attain sizes smaller than $\eta$ was studied in
\cite{bermar}. It was shown that   the
mean energy behaves as $1/ \eta^{\alpha-\beta}$ when $\eta\to0$,
where $\alpha$ denotes the Malthusian exponent of the
fragmentation process and where $\alpha>\beta$ in physically
reasonable cases. Therefore, the performances of two individual
fragmentation processes are asymptotically comparable by means of
the quantities $\alpha-\beta$ and $\widehat{\alpha}-\widehat{\beta}$,
where $\widehat{\alpha}>\widehat{\beta}$ are the parameters associated
with a second fragmentation process.

We shall  formulate problem (*) in  mathematical terms adopting
the same mean energy point of view as in \cite{bermar}. First, we
will explicitly compute the objective function, which we express
in terms of the Levy and renewal measures associated with the
``tagged fragment'' of each of the two fragmentation processes
(see \cite{ber2}). Then, our  goal   will be to  study  a
preliminary question related to (*), which is weaker but still
relevant  for the mining industry:

\smallskip

(**) \ when is the above described ``two-step'' procedure
efficient in terms of mean energy, compared to the ``one-step''
procedures where only the first or only the second fragmentation
mechanisms reduce a unit size fragment to fragments not larger
that $\eta_0$?

\smallskip

We shall address this question in  asymptotic regimes, namely for
$\eta$ and $\eta_0$ going together either to $0$ or to $1$. In
both cases, we will give explicit estimates in terms of $\eta$ for
the efficiency gain or loss of using the two-step procedure.

As we shall see, if  $\alpha,\beta$, $\widehat{\alpha}$ and
$\widehat{\beta}$ are different, for any values of $\eta/\eta_0\in
(0,1)$ the relations between those four parameters determine  the
relative efficiency between the first, the second, and the
two-step fragmentation procedures if $\eta$ is sufficiently small.
In particular, when  $\alpha>\widehat{\alpha}$ and $\beta
>\widehat{\beta}$ the answer to question (**) is affirmative for
$\eta$ sufficiently small, so that the solution to problem (*) is
in general  non trivial.

We shall carry out a similar analysis for large (that is, close to
unit-size) thresholds.  In order to quantify the comparative
efficiency of the two-step procedure, we shall make an additional
hypothesis of  regular variation at $\infty$ of the Levy exponents
of the tagged fragment  processes. This will be transparently
interpreted in terms of the infinitesimal average energy required
by each of the fragmentation processes to break arbitrarily close
to unit-size fragments. We will show that at least for small
values of $\log \eta_0 /\log \eta$ and variation indexes in
$(0,\frac{1}{2}]$ for both fragmentation processes, the relative
infinitesimal efficiency of the two fragmentation processes
determines the comparative efficiency of the three alternative
fragmentation procedures if $\eta$ is sufficiently close to $1$.

We point out that the relevant parameters involved in our analysis
could in principle be statistically estimated. A first concrete
step in that direction has been made by Hoffmann and Krell
\cite{kre} who asymptotically estimate the Levy measure of the
tagged fragment from the observations of the sizes fragments at
the first time they become smaller than $\eta\to 0$. Although this
is in general not enough to recover the characteristics of the
fragmentation process, it provides all the relevant parameters we
need which are not observable by other means.

The remainder of this paper is organized as follows. In Section 2
we recall the construction of homogeneous fragmentation processes
in terms of Poisson point processes, we describe our model of the
two-step fragmentation procedure and compute its average energy
using first passage laws for subordinators. In Section 3 we recall
some results on renewal theory for subordinators and use them to
study the small thresholds asymptotics of our problem in Theorems
\ref{2frvs1} and \ref{2frvs2}, where the two-step procedure is
respectively compared with the first and the second fragmentation
processes. The comparative efficiency of the three alternatives
according to the values of $\alpha,\beta$, $\widehat{\alpha}$ and
$\widehat{\beta}$ is summarized in Corollary \ref{summary}.  In
Section 4 we introduce the idea of relative ``infinitesimal
efficiency'' of two fragmentation procedures. We relate it to a
regular variation assumption at infinity for  the Levy exponent of
tagged fragment, and use it to analyze the comparative efficiency
of the two-step fragmentation procedure for close to unit-size
fragments, using Dynkin-Lamperti asymptotics for subordinators at
first passage.

\section{The model}
\subsection{The fragmentation process}

We shall model the fragmentation mechanisms as a homogeneous
fragmentation processes, as introduced in \cite{ber}.
This is a homogeneous Markov process $\X=(X(t,\x) :t\geq 0)$
taking values in $$\S:=\left\{\s= (s_{1},s_{2},...) \ : \
s_{1}\geq s_{2}\geq ...\geq 0\ , \sum_{i=1}^{\infty} s_{i} \leq 1
\right\}\,,
$$
which satisfies the two fundamental properties of homogeneity and
branching.
The parameter $\x=(x_{1},x_{2},...)$ is an element of
$\S$ standing for the initial condition:
$X(0,\mathbf{x} )=\x$ a.s.. In the case $\x=(1,0,\dots)$ we
simply write
$X(t)=X(t,\x)$, $t\geq 0$.

\medskip

We observe that homogeneous fragmentation processes are
self-similar fragmentation processes with
zero index of self-similarity  (see \cite{ber}). Since
self-similar fragmentation processes with different indexes are
related by a family of random time-changes (depending on
fragments), there is no loss of generality in working here in the
homogeneous case as the quantities we study are only
size-dependent (see also \cite{bermar}).

\medskip

We  assume that no creation of mass occurs. It is known that in
this case, the process $\X$ is entirely characterized by an
erosion coefficient $c\geq 0$ and a dislocation measure $\nu$,
which is a measure on $\S$ satisfying the conditions
\begin{equation}
\label{mesuredelevy}
\nu (\{1,0,0,...\})=0 \, \hbox{ and } \, \int_{\S}
(1-s_{1})\nu(d\s)<\infty\,.
\end{equation}
Moreover, we suppose that we are in the dissipative case
$\sum_{i=1}^{\infty}s_{i}\leq 1 \hbox{ a.s.}$, and we assume
absence of erosion: $c=0$.

\medskip

Let us recall the construction of a homogeneous fragmentation
process in this setting, in terms the atoms of a Poisson point
process (see \cite{be2}). Let $ \nu$ be a dislocation measure
fulfilling conditions \eqref{mesuredelevy}. Let $\K = \left(
(\Delta(t) , k(t)): t \geq 0 \right)$ be a Poisson point process
with values in $\S\times \NN$, and with intensity measure
$\nu\otimes\sharp$, where $\sharp$ is the counting measure on
$\NN$. As in \cite{be2}, we can construct a unique $\S$-valued
process $\X= (X(t,\x): t\geq 0)$ started from $\x$ with paths that
jump only at instants $t\geq 0$ at which a point $
(\Delta(t)=(\Delta_{1},\Delta_{2},....),k (t))$ occurs. Plainly,
$X(t,\x)$ is obtained by replacing the $k(t)$-term $ X(t-,\x)$ by
the decreasing rearrangement of the sequence
$X_{1}(t-,\mathbf{x}),...,X_{k-1}(t-,\x),X_{k}(t-,\x)
\Delta_{1},X_{k}(t-,\x)\Delta_{2},...,X_{k+1}(t-,\x),... $.

\medskip

Define
$$\underline{p}:=\inf\left\{p\in\mathbb{R}:\
\int_{\S}\sum_{j=2}^{\infty}s_{j}^{p} \nu (d\mathbf{s})<\infty \right\}
$$
and for every $q\in (\underline{p},\infty)$ consider,
\begin{equation}
\label{kappa}
\kappa (q) :=
\int_{\S}\left (1-\sum_{j=1}^{\infty}s_{j}^{q}\right)
\nu (d\s)\, .
\end{equation}

In the sequel, we assume the Malthusian hypothesis:
$\exists \, ! \, \alpha \geq \underline{p}$ such that
$\kappa (\alpha )=0$ which is called the Malthusian exponent.

\medskip

A key tool in fragmentation theory is
the tagged fragment associated with $\X$. For
the precise definition, we refer the reader to \cite{ber2}.
The tagged fragment is a process defined by
$$
\chi (t):=X_{J(t)}(t)
$$
where $J(t)$ is a random integer such
that, conditioned on $X(t)$, $\mathbb{P}(J(t)=i|X(t))=X_{i}(t)$ for all
$i\geq 1$, and $\mathbb{P}(J(t)=0|X(t))=1-\sum_{i=1}^{\infty}X_{i}(t)$.

\medskip

Is is shown by Bertoin (Theorem 3 in \cite{ber2}) that the process
$$
\xi_{t}=-\log \chi (t)
$$
is a subordinator. Moreover, its Laplace exponent $\phi$ is given by
$$\phi(q):=\kappa (q+1)$$ for $q>\underline{p}-1$.
Since $\phi(\alpha-1)=0$, the process $e^{(1-\alpha)\xi (t)}$ is a
nonnegative martingale, and we can then define a probability measure
$\widetilde{\mathbb{P}}$ on the path space by
\begin{equation}
\label{eq1111}
d\widetilde{\mathbb{P}}{\big|}_{\mathcal{F}_{t}}=e^{(1-\alpha)\xi(t)}
d\mathbb{P}{\big|}_{\mathcal{F}_{t}} \,,
\end{equation}
where $\left(\mathcal{F}_{t} : t\geq 0 \right)$ denotes the
natural filtration of $\xi$. It is well known that under this
``tilted'' law, $\xi$ is a subordinator with Laplace exponent
\begin{equation}
\label{phitilde}
\widetilde{\phi}(q)=\phi(q+\alpha-1).
\end{equation}

We will respectively denote by ${\Pi}$ and
${U}$ the L\'evy measure and the renewal measure of
$\xi_{t}$ under ${\widetilde{\mathbb{P}}}$ (see e.g. \cite{be1}).

\medskip

For $\eta\in (0,1]$ we denote by
$$
T_{\eta}:=\inf\{ t\geq 0 : \xi_{t}>\log(1/\eta)\}
$$
the first time that the size of the tagged fragment is smaller than
$\eta$.

\subsection{ The fragmentation energy }
Following \cite{bermar}, we shall assume that the
energy needed to split a fragment of size $x\in [0,1]$ into a sequence
$x_1\geq x_2\geq \dots$ is given by the formula
$$x^{\beta}\varphi\left(\frac{x_1}{x},\frac{x_2}{x},\dots\right),$$ where
$\beta>0$ is a fixed constant and
 $\varphi :\mathcal{S}\rightarrow \mathbb{R}$ is a measurable ``cost
function'' such that $\varphi ((1,0,...))=0$.

\medskip

We are interested in the total energy ${\mE}^{(\x)} (\eta)$ used
in splitting the initial fragment of size $x$
until each of them has reached, for the first
time, a size that is smaller than $\eta$. This
quantity is given by
$$
{\mE}^{(\mathbf{x})}(\eta)=\sum_{t\geq 0}
\nbOne_{X_{k(t)}(t_{-}, \x)\geq \eta} X_{k(t)}^{\beta}
(t_{-},\mathbf{x})\varphi (\Delta (t)).
$$
We shall simply write
$$
{\mE}(\eta):={\mE}^{(1,0,\dots)}(\eta)\,.
$$

The following consequence of the homogeneity property will be
useful.

\begin{lm}
\label{ssproperty}
Let $\mathbf{x}=(x_{1},x_{2},...)\in
\mathcal{S^{\downarrow}}$ and $\eta \in [0,1]$. We have
\begin{equation}
\label{eqenergie} {\mE}^{(\mathbf{x})}(\eta^{})
\stackrel{(law)}{=} \sum_{i} \nbOne_{x_{i}\ge \eta^{}}
x_{i}^{\beta^{}}{\mE}_{i}(\eta^{}/x_{i}),
\end{equation}
where for each $i\geq 1$, ${\mE}_{i}(\cdot)$ is the energy of a
fragmentation process $\X^{(i)}$ issued from $(1,0,\dots)$ with the same
characteristics as $\X$, and the copies $(\X^{(i)}: i\geq 1)$ are
independent.
\end{lm}

\begin{dem}
Let $\left((\Delta_{i}(t),k_{i}(t)): t\geq 0 \right)$, $i\geq 1$,
be i.i.d. Poisson point processes with intensity measure $\nu\otimes
\sharp$. Denote by  ${\overline {\X}}^{(x_{i})}$, $i\geq 1$,  the
sequence of
independent homogeneous fragmentation processes constructed from
the latter processes, respectively  starting from
$(x_{i},0,\cdots)$.  From the branching property of $\X$, we have
the identity
$$
{\mE}^{(\x)}(\eta^{})\stackrel{(law)}{=}\sum_{i}\sum_{t\geq 0}
\nbOne_{x_{i}\geq\eta^{}}\, \nbOne_{{\overline X}_{k_{i}(t)}^{(x_{i})}
(t_{-})\geq \eta^{}}
({\overline X}_{k_{i}(t)}^{(x_{i})})^{\beta^{}}(t_{-})\,
\varphi^{} (\Delta_{i} (t)).
$$
Denoting now by $\left( (\Delta^{(i)}(t),k^{(i)}(t)): t\geq 0
\right)$, $i\geq 1$, the family of i.i.d. Poisson point processes
associated with the process $\X^{(i)}$, we get by  homogeneity
that
$$
{\mE}^{(\x)}(\eta^{})\stackrel{(law)}{=}\sum_{i}\sum_{t\geq 0}
\nbOne_{x_{i}\geq\eta^{}} \,
\nbOne_{x_{i}X_{k^{(i)}(t)}^{(i)}(t_{-})\geq \eta^{}}
x_{i}^{\beta^{}} (X_{k^{(i)}(t)}^{(i)})^{\beta^{}}(t_{-}) \,
\varphi^{} (\Delta^{(i)} (t)),
$$
and the statement follows.
\end{dem}

\subsection{ The energy of a two-step  fragmentation procedure }

To formulate our problem, we introduce a second Poisson point
process ${\widehat \K} = (({\widehat \Delta}(t), {\widehat k}(t)
), t \geq 0 )$ with values in $\S\times \NN$, and with intensity
measure $\hnu\otimes\sharp$, where $\hnu$ is a dislocation measure
satisfying the same type of assumptions as $\nu$. We can then
simultaneously define a family of fragmentation processes
$\hX=({\widehat X}(t,\x): t\geq 0)$ indexed by the initial
condition $\x=(x_{1},x_{2},...)$. We denote by $\ba$ the Malthus
coefficient of $\hnu$. The energy used in the second fragmentation
process is assumed to take the same form as for the first, in
terms of (possibly different) parameters $\bb$ and $\hvi$.

\medskip

We assume that $\K$ and ${\widehat \K}$ are independent, so the
families of fragmentation processes $\X$ and  $\hX$ are
independent, and they are called respectively the first and the
second fragmentation processes.

\medskip

In the sequel we assume that the first fragmentation process
$\X$ is issued from the unitary
fragment $(1,0,\dots)$. Let $1\geq \eta \geq \eta_0>0$.
We let each mass fragment evolve in the first fragmentation
process until the instant it first becomes smaller than
$\eta$. Then it immediately enters the second
fragmentation process $\hX$, and then evolves until it first becomes
smaller than $\eta_0$.

\medskip

For each $\eta\in (0,1]$ let $\x^{\eta}\in \mathcal{S^{\downarrow}}$
be the mass partition given by the ``output''
of $\X$ when each of the fragments reaches for the first time a size
smaller than $\eta$. More precisely, each fragment is
``frozen'' at that time, while other (larger than $\eta$) fragments
continue their independent evolutions. We write
\begin{equation}
\label{ouput1}
\x^{\eta}=(x_1^{\eta},x_2^{\eta},\cdots)
\end{equation}
for the decreasing rearrangement of the (random) frozen sizes of
fragments when exiting the first fragmentation process.
By the homogeneity and
branching properties, if $\E(\eta ,\eta_0)$ denotes the total
energy spent in reducing the unit-size fragment by these
procedure, we have the identity
\begin{equation}
\label{energieloi}
\E(\eta,\eta_0):\stackrel{(law)}{=} {\mE}(\eta)+
{\mmE}^{(\x_{\eta})}(\eta_0)\,,
\end{equation}
where ${\mmE}^{(\x)}(\cdot)$ is the energy of a copy of the
second fragmentation process $\hX$ starting from $\x$,
independent of the first fragmentation process.

\medskip

\begin{rk}
Notice that $\E(1,\eta_0)$ is the energy required to initially
dislocate the  unit mass  with the first fragmentation process,
and then use the second fragmentation process to continue breaking
its fragments if their sizes are larger or equal to $\eta_0$ (the
other ones immediately exit from the system). We will denote
$\E(1^+,\eta_0)={\mmE}(\eta_0)$ the total energy required when
only the second fragmentation process is used from the beginning.

For the quantity $\E(\eta_0,\eta_0)={\mE}(\eta_0)$ no confusion
arises: it corresponds to the case when the first fragmentation
process is used during the whole procedure.
\end{rk}

Our goal now is to compute the expectation of $\E(\eta ,\eta_0)$.

\medskip

The notation ${\widehat \xi}$, ${\widehat T}_{\eta}$,
${\widehat{\Pi}}$, ${\widehat {U}}$ and so on, will be used for
the analogous objects associated with the fragmentation process
$\hX$.

So far the notation ${\mathbb{P}}$ has been used to denote the law
of $\xi$. In all the sequel, we keep the same notation
${\mathbb{P}}$  to denote the product law of independent copies of
the processes $\xi$ and $\widehat{\xi}$ in the product path space.
Extending accordingly the  definition in (\ref{eq1111}), we will
also denote by $\widetilde{\mathbb{P}}$  the product measure the
first marginal  of which is given by
$d\widetilde{\mathbb{P}}{\big|}_{\mathcal{F}_{t}}=e^{(1-\alpha)\xi(t)}
d\mathbb{P}{\big|}_{\mathcal{F}_{t}}$ and the second one given by
$d\widetilde{\mathbb{P}}{\big|}_{{\widehat {\mathcal{F}}}_{t}}=
e^{(1-{\widehat \alpha}){\widehat \xi}(t)}
d\mathbb{P}{\big|}_{{\widehat {\mathcal{F}}}_{t}}$. Here
$\left(\mathcal{F}_{t} : t\geq 0 \right)$ and $\left({\widehat
{\mathcal{F}}}_{t} : t\geq 0 \right)$ are the natural filtrations
of $\xi$ and $\widehat \xi$ respectively.

\medskip

We shall assume throughout that
the following integrability condition holds:
\begin{equation}
\label{inte1}
\varphi\in L^1(\nu)\; \hbox{ and } \; \hvi\in L^1(\hnu)\,.
\end{equation}
In this case we define
$$
C= \int_{ S } \varphi(\s) \nu( d \s) \; \hbox{ and } \;
\cc= \int_{ S } \hvi(\s) \hnu( d \s)\,.
$$
Let us introduce the functions
$$
\Psi(x)=C \, \int_{0}^{x} e^{(\alpha-\beta) y} {U}
(dy)\,,\;\;\;
\hPs(x)=
\cc \,\int_{0}^{x} e^{(\ba-\bb) y} \widehat{U} (dy)
\,,\; x\geq 0\, .
$$
To simplify the notation we will put
$$
\forall\,  a>0 : \;\, \ell(a):=\log (1/a)\,.
$$

We have the elements to compute the expected energy requirement in
the two step fragmentation procedure.

\begin{lm}
\label{expen} Assume that the integrability condition
(\ref{inte1}) is satisfied. Let $\eta_0\in (0,1)$. Then, we have
for $\eta_0<\eta<1$ that
$$
\begin{disarray}{rcl}
 \mathbb{E}(\E(\eta ,\eta_0))&=&
 C \,\int_{0}^{ \ell(\eta)} e^{(\alpha-\beta) y} \, {U}
(dy)
\\
&& + \cc \, \int_{0}^{ \ell(\eta)}
\int_{ \ell(\eta)-y}^{\ell(\eta_0)-y}
e^{(\alpha-\bb)(z+y)} \left[\int_{0}^{\ell(\eta_0)-(z+y)}
e^{(\ba-\bb)x}
\widehat{U}(dx)\right] {\Pi}(dz){U}(dy)
 \\
&=&\Psi( \ell(\eta))+ \widetilde {\mathbb{E}}
\left(\nbOne_{\xi_{T_{\eta}}< \ell(\eta_0)} e^{(\alpha-\bb)
\,\xi_{T_{\eta}}} \, \hPs( \ell(\eta_0)-\xi_{T_{\eta}})\right),
\\
\end{disarray}
$$
and
$$
\mathbb{E}(\E(\eta_0 ,\eta_0))= \mathbb{E}({\mE}(\eta_0))=
\Psi(\ell(\eta_0)), \quad \mathbb{E}(\E(1^+ ,\eta_0))=
\mathbb{E}({\mmE}(\eta_0))=\hPs( \ell(\eta_0)).
$$
When the renewal measures $U(dx)$ has no atom at $0$ one has
$$
\mathbb{E}(\E(1^+ ,\eta_0))= \mathbb{E}(\E(1 ,\eta_0))\,.
$$
\end{lm}

\begin{dem} The proof is an extension of arguments given in
\cite{bermar} corresponding to the case ``$\eta=1^+$'' or
$\eta=\eta_0$ and which we repeat here for convenience.  By the
compensation formula for the Poisson point process $(\Delta (u),
k(u))$ associated with the first fragmentation process $\X$, we
get that for $\eta_0\in (0,1]$,
\begin{equation*}
\label{eq1}
\mathbb{E}({\mE}(\eta_0))
=\mathbb{E}\left(\int_{0}^{\infty} \nbOne_{\chi(t)>\eta_0}
(\chi (t))^{\beta-1} dt\right) \int_{ S } \varphi (\mathbf{s})
\nu ( d \s)
=C ~\mathbb{E}\left(\int_{0}^{\infty}
\nbOne_{\xi_{t}<\ell(\eta_0) }
e^{(1-\beta)\xi_{t}} dt\right)\,.
\end{equation*}
Thus
\begin{equation}
\label{aven}
\mathbb{E}({\mE}(\eta_0))
=C\widetilde{\mathbb{E}}\left(\int_{0}^{\infty}
\nbOne_{\xi_{t} <\ell(\eta_0) }
e^{(\alpha-\beta)\xi_{t}} dt\right)
=C\int_{0}^{\ell(\eta_0)} e^{(\alpha-\beta) y} {U} (dy)
=\Psi(\ell(\eta_0)).
\end{equation}
Similarly,
$$
\begin{disarray}{rcl}\mathbb{E}({\mmE}(\eta_0))
&=&\cc \int_{0}^{\ell(\eta_0)} e^{(\ba-\bb) y} \widehat{U}
(dy)=\hPs(\ell(\eta_0)).
\end{disarray}
$$
The above identity also implies that $\mathbb{E}(\E(1 ,\eta_0))=
\cc \int_{0^+}^{\ell(\eta_0)} e^{(\ba-\bb) y} \widehat{U}(dy)=
\mathbb{E}({\mmE}(\eta_0))$ when $U$ has no atom at $0$.

\medskip

The statement is thus proved for the cases ``$\eta=1^+$'' and
$\eta=\eta_0$. For the general case, we use Lemma \ref{ssproperty}
to get
$$
\begin{disarray}{rcl}
\mathbb{E}({\mmE}^{(\mathbf{x}_{\eta})}(\eta_0))&=&
\mathbb{E}\left(\sum_{i} \nbOne_{x_{\eta,i}>\eta_0}
x_{\eta,i}^{\bb}{\mmE}_{ i}(\eta_0/x_{\eta,i})\right)\\&=&
\mathbb{E}\left(\sum_{i} \nbOne_{x_{\eta,i}>\eta_0}
x_{\eta,i}^{\bb}\mathbb{E}
({\mmE}_{i}(\eta_0/x_{\eta,i}) \, | \, x_{\eta ,i})\right)\\
&=& \mathbb{E}
\left( \nbOne_{\chi(T_{\eta})>\eta_0}(\chi (T_{\eta}))^{\bb-1}
{\widehat {\mathbb{E}}}({\mmE}(\eta_0/y))\vert_{y= \chi (T_{\eta})} \right),
\end{disarray}
$$
where
${\mmE}(\cdot)$ is the energy of a copy of the second fragmentation
process, starting from the unit mass, and which is independent of the
first one, and ${\mmE}_{i}(\cdot)$ are independent copies of
${\mmE}(\cdot)$.
Then, since $\chi (t)=e^{-\xi_t}$, we have,
$$
\begin{disarray}{rcl}
\mathbb{E}({\mmE}^{(\mathbf{x}_{\eta})}(\eta_0))&=&
\widetilde {\mathbb{E}} \left( \nbOne_{\xi_{T_{\eta}}<
\ell(\eta_0)}e^{(\alpha-\bb)\xi_{T_{\eta}}}
\widetilde{\mathbb{E}}({\mmE}(\eta_0e^{z} ))\vert_{
z=\xi_{T_{\eta}}}\right) \\ &=& \widetilde {\mathbb{E}}
\left(
\nbOne_{\xi_{T_{\eta}}<\ell(\eta_0)}e^{(\alpha-\bb)
\xi_{T_{\eta}}}\hPs(\ell(\eta_0)-\xi_{T_{\eta}})\right)
.\\
\end{disarray}
$$
According to Lemma 1.10 of \cite{ber1999} the
distribution of $\xi_{T_{\eta}}$ under
$\widetilde{\mathbb{P}}$ is given by
$$
\widetilde{\mathbb{P}}(\xi_{T_{\eta}}\in dz)= \int_{0}^{\ell(\eta)}
\nbOne_{ \ell(\eta )< z} {\Pi} (dz-y){U} (dy) .
$$
Therefore,
$$
\mathbb{E}({\mmE}^{(\x_{\eta})}(\eta_0))=
\int_{0}^{ \ell(\eta)}
\left[\int_{\ell(\eta)-y}^{ \ell(\eta_0)-y}e^{(\alpha-\bb)(z+y)}
\hPs \left( \ell(\eta)-(z+y)\right){\Pi}(dz)\right]{U}(dy).
$$
By bringing the pieces together and by using the identity
\eqref{energieloi} we get the result.
\end{dem}

\medskip

In analogy with (\ref{ouput1}), we introduce the notation
\begin{equation}
\label{ouput2}
{\hx}^{\eta}=({{\widehat x}_1}^{\eta},{{\widehat x}_2}^{\eta},\cdots)
\end{equation}
for the decreasing rearrangement of the frozen sizes of fragments
smaller than $\eta$, that exit the second
fragmentation process started from the unit mass. The following
decompositions of the total energy will be useful in the sequel:

\begin{rk}
\label{decomp}
For $1\geq \eta \geq \eta_0 >0$ we have
$$
{\mmE}(\eta_0)= {\mmE}(\eta)+
{\mmE}^{({\hx}^{\eta})}(\eta_0),
$$
whence,
\begin{equation*}
\begin{split}
\E(\eta ,\eta_0)-\E(1^+,\eta_0)= & {\mE}(\eta) -{\mmE}(\eta)
+ {\mmE}^{(\x^{\eta})}(\eta_0) -
{\mmE}^{({\hx}^{\eta})}(\eta_0)\,.
\end{split}
\end{equation*}
From this relation and by similar computations as in Lemma
\ref{expen}, we can write
\begin{equation*}
\begin{split}
\mathbb{E}(\E(\eta ,\eta_0)-\E(1^+,\eta_0))
&= \Psi(\ell(\eta)) - \hPs(\ell(\eta)) \\
&\quad + {\widetilde {\mathbb{E}}}
\left( \nbOne_{\xi_{T_{\eta}}< \ell(\eta_0)}
e^{(\alpha-\bb)\xi_{T_{\eta}}}
\hPs(\ell(\eta_0)-\xi_{T_{\eta}})\right) \\
& \quad -{\widetilde{\mathbb{E}}}
\left(\nbOne_{ {\widehat \xi}_{ {\widehat T}_{\eta} }<\ell(\eta_0) }
e^{(\ba-\bb){\widehat \xi}_{ {\widehat T}_{\eta} } }
\hPs(\ell(\eta_0)-{\widehat \xi}_{ {\widehat T}_{\eta} } )
\right)
.\\
\end{split}
\end{equation*}
Observe that when $U$ has no atom at $0$, one can replace
$\E(1^+,\eta_0)$ by $\E(1,\eta_0)$ on the left
hand side of the formula.

Similarly, we have
\begin{equation*}
\begin{split}
\mathbb{E}(\E(\eta ,\eta_0)- &\E(\eta_0,\eta_0))=
{\mmE}^{(\x^{\eta})}(\eta_0) -
{\mE}^{({\x}^{\eta})}(\eta_0)\\ = &
{\widetilde {\mathbb{E}}}
\bigg( \nbOne_{\xi_{T_{\eta}}<
\ell(\eta_0)}e^{(\alpha-\bb)\xi_{T_{\eta}}}
\hPs(\ell(\eta_0)-\xi_{T_{\eta}})  - \nbOne_{\xi_{T_{\eta}}<
\ell(\eta_0)}e^{(\alpha-\beta)\xi_{T_{\eta}}}
\Psi(\ell(\eta_0)-\xi_{T_{\eta}})\bigg).\\
\end{split}
\end{equation*}
\end{rk}

\section{Small thresholds}

In this  section, we consider the  total energy
$\mathbb{E}(\E(\eta ,\eta_0))$ when  $\eta_0$ and $\eta$ go to $0$
in a suitable joint asymptotics. Our goal is to compare it with
the mean energy required for reducing the unit fragment to
fragments smaller than $\eta_0$ using  only the first or only the
second fragmentation processed.  We shall assume that the
quantities

$$
m(\alpha):=\int_{\S}\sum_{i=1}^{\infty}s_n^{\alpha}
\log\left(\frac{1}{s_n^{\alpha}}\right)\nu(d\s)\, \mbox{ and } \,
{\widehat m}(\ba):=\int_{\S}\sum_{i=1}^{\infty}s_n^{\alpha}
\log\left(\frac{1}{s_n^{\ba}}\right)\hnu(d\s)
$$
are finite. Moreover, we impose the conditions
$$
\beta<\alpha\mbox{  and }\bb<\ba.
$$
The latter assumption is  physically reasonable, since the energy
$\Psi(\infty)$ (respectively $\hPs(\infty)$) required in order
that all fragments vanish in the first (respectively second)
fragmentation processes is otherwise finite (see Remark 1 in
\cite{bermar}).

\medskip

The following asymptotic result on the mean energy of a single
fragmentation processes is based on the renewal Theorem for
subordinators (Bertoin et al. \cite{Be99renew}). Its proof is
simply adapted from that of Lemma 4 in \cite{bermar}, see also
Theorem 1 therein.

\begin{lm}
\label{ensmallfrag} Under the previous assumptions, we have
$$\lim_{\eta\to 0} \eta^{\alpha-\beta}
\mathbb{E}({\mE}(\eta))=\frac{C}{(\alpha-\beta)m(\alpha)} \,
\mbox{ and } \, \lim_{\eta\to 0} \eta^{\ba-\bb}
\mathbb{E}({\mmE}(\eta))=\frac{\cc}{(\ba-\bb){\widehat m}(\ba)}.$$
\end{lm}

By the renewal Theorem for subordinators we also have as $\eta\to
0^+$ that
$$
\widetilde{\mathbb{P}} \left(\xi_{T_{\eta}}-\ell(\eta)\in
du\right) \to M(du):=\frac{1}{m(\alpha)}\int_{\mathbb{R}^+}
{\Pi}(y+du)dy,
$$
and
$$
\widetilde{\mathbb{P}} \left({\widehat \xi}_{{\widehat T}_{\eta}}-
\ell(\eta)\in du\right) \to {\widehat M}(du):= \frac{1}{{\widehat
m}(\ba)}\int_{\mathbb{R}^+} {\widehat \Pi}(y+du)dy
$$
in the weak sense.  Let us define, for $\lambda>0$  a fixed
parameter, the finite and strictly positive constants

$$
F_{\lambda}:=m(\alpha)\int_0^{\lambda}e^{(\alpha-\bb)u}
\hPs(\lambda-u) M(du),$$
$$D_{\lambda}:=m(\alpha)\int_0^{\lambda}e^{(\alpha-\beta)u}
\Psi(\lambda-u) M(du)$$
$$
{\widehat D}_{\lambda}:={\widehat
m}(\ba)\int_0^{\lambda}e^{(\ba-\bb)u} \hPs(\lambda-u) {\widehat
M}(du).
$$

We fix in the sequel the parameter $\lambda>0$.  With these
elements, we are in position to explicitly study the (comparative)
behavior of the total energy for small thresholds $\eta$ and
$\eta_0$, when these are bond by the relation
$$
\eta_0=\eta e^{-\lambda}.
$$

\begin{theo}
\label{2frvs1} (Two-step procedure versus first fragmentation
only)

Assume that the renewal measure $U(dx)$ has no atom at $0$.
 For any $\lambda>0$, the following hold:

$(a)$ If $\bb>\beta$, then $\forall \, \varepsilon\in
(0,\frac{\alpha-\beta}{C}D_{\lambda})\;$ $\exists \,
\eta_{\varepsilon}\in (0,1)$ such that
$$
\forall \, \eta\leq \eta_{\varepsilon}: \;\;\; \mathbb{E}(\E(\eta
,\eta e^{-\lambda}))\leq
\left(\varepsilon-\frac{\alpha-\beta}{C}D_{\lambda}\right)
\mathbb{E}({\mE}(\eta))+\mathbb{E}({\mE}(\eta e^{-\lambda}))<
\mathbb{E}({\mE}(\eta e^{-\lambda})).
$$
\noindent $(b)$ If $\bb<\beta$, then $\forall \, M>0\;$ $\exists
\, \eta_{M}\in (0,1)$ such that
$$
\forall \, \eta\leq \eta_{M}: \;\;\; \mathbb{E}(\E(\eta ,\eta
e^{-\lambda}))\geq M \mathbb{E}({\mE}(\eta))
+\mathbb{E}({\mE}(\eta e^{-\lambda}))> \mathbb{E}({\mE}(\eta
e^{-\lambda})).
$$
\noindent $(c)$ If $\bb=\beta$, then $\forall \, \varepsilon\in
(0,1)$ $\exists \, \eta_{\varepsilon}\in (0,1)$ such that $\forall
\, \eta\leq \eta_{\varepsilon}$:
 \begin{multline*}
\left(-\varepsilon+\frac{\alpha-\beta}{C}(F_{\lambda}-D_{\lambda})\right)
\mathbb{E}({\mE}(\eta))+\mathbb{E}({\mE}(\eta e^{-\lambda})) \\
\leq \mathbb{E}(\E(\eta ,\eta e^{-\lambda}))\leq
\left(\varepsilon+\frac{\alpha-\beta}{C}(F_{\lambda}-D_{\lambda})\right)
\mathbb{E}({\mE}(\eta))+\mathbb{E}({\mE}(\eta e^{-\lambda}))\,.
\end{multline*}

In all cases, one can replace $\mathbb{E}({\mE}(\eta))$ by
$C\left[ \eta^{\alpha-\beta} m(\alpha)(\alpha-\beta)\right]^{-1}$.
\end{theo}

\begin{dem}
All parts are obtained by taking limit when $\eta\to 0$ in the
identity
\begin{equation*}
\begin{split}
\eta^{\alpha-\beta}\left(\mathbb{E}(\E(\eta ,\eta
e^{-\lambda}))\!-\!\mathbb{E}(\E(\eta e^{-\lambda} ,\eta
e^{-\lambda}))\right)= & \widetilde {\mathbb{E}} \bigg(
e^{(\alpha-\bb)(\xi_{T_{\eta}}-\ell(\eta))}
\hPs(\lambda\!-\!(\xi_{T_{\eta}}\!-\!\ell(\eta))
\nbOne_{\xi_{T_{\eta}}-\ell(\eta)<\lambda}\bigg)\eta^{\bb-\beta}  \\
& -\widetilde {\mathbb{E}} \bigg(
e^{(\alpha-\beta)(\xi_{T_{\eta}}-\ell(\eta))}
\Psi(\lambda\!-\!(\xi_{T_{\eta}}\!-\!\ell(\eta))
\nbOne_{\xi_{T_{\eta}}-\ell(\eta)<\lambda }\bigg),
\end{split}
\end{equation*}
which follows from Remark \ref{decomp}, and then using Lemma
\ref{ensmallfrag} and the previously mentioned weak convergence
result for $\widetilde{\mathbb{P}} \left({\widehat \xi}_{{\widehat
T}_{\eta}} -\ell(\eta)\in dy\right)$ (notice that the limit is
absolutely continuous).
\end{dem}

\medskip

\begin{theo} \label{2frvs2} (Two-step procedure versus second fragmentation only)

Assume that $U(dx)$ and $\widehat{U}(dx)$ have no atom at $0$.
For any $\lambda>0$, the following hold:

$(a)$ If $\ba>\alpha$, then $\forall \, \varepsilon\in
(0,\frac{\ba-\bb}{\cc} {\cd}_{\lambda})\;$ $\exists \,
\eta_{\varepsilon}\in (0,1)$ such that
$$
\forall \, \eta\leq \eta_{\varepsilon} : \;\;\; \mathbb{E}(\E(\eta
,\eta e^{-\lambda}))\leq
\left(\varepsilon-\frac{\ba-\bb}{\cc}{\cd}_{\lambda}\right)
\mathbb{E}({\mmE}(\eta))+ \mathbb{E}({\mmE}(\eta e^{-\lambda}))<
\mathbb{E}({\mmE}(\eta e^{-\lambda}))\,.
$$

\noindent $(b)$ If $\ba<\alpha$, then $\forall \, M>0 \;$ $\exists
\, \eta_{M}\in (0,1)$ such that
$$
\forall \eta\leq \eta_{M}: \;\;\; \mathbb{E}(\E(\eta ,\eta
e^{-\lambda}))\geq M \mathbb{E}({\mmE}(\eta))
+\mathbb{E}({\mmE}(\eta e^{-\lambda}))> \mathbb{E}({\mmE}(\eta
e^{-\lambda})).
$$

\noindent $(c)$ If $\ba=\alpha$, then $\forall \, \varepsilon\in
(0,1)\;$ $\exists \, \eta_{\varepsilon}$ such that $\forall \,
\eta\leq \eta_{\varepsilon}$,
\begin{multline*}
\left(-\varepsilon+\frac{\ba-\bb}{\cc}(F_{\lambda}-{\cd}_{\lambda})\right)
\mathbb{E}({\mmE}(\eta))+\mathbb{E}({\mmE}(\eta e^{-\lambda}))
\\ \leq \mathbb{E}(\E(\eta ,\eta e^{-\lambda}))\leq
\left(\varepsilon+\frac{\ba-\bb}{\cc}(F_{\lambda}-{\cd}_{\lambda})\right)
\mathbb{E}({\mmE}(\eta))+\mathbb{E}({\mmE}(\eta e^{-\lambda}))\,.
\end{multline*}

In all cases, one can replace $\mathbb{E}({\mmE}(\eta))$ by
$\cc\left[\eta^{\ba-\bb}m(\ba)(\ba-\bb)\right]^{-1}$.
 \end{theo}

\begin{dem}
The proof is similar to previous one, noting that
\begin{equation*}
\begin{split}
\eta^{\ba-\bb}\left(\mathbb{E}(\E(\eta ,\eta e^{-\lambda}))
\!-\!\mathbb{E}(\E(1 ,\eta e^{-\lambda}))\right)= & \widetilde
{\mathbb{E}} \bigg( e^{(\alpha-\bb)(\xi_{T_{\eta}}-\ell(\eta))}
\hPs(\lambda\!-\!(\xi_{T_{\eta}}\!-\!\ell(\eta))
\nbOne_{\xi_{T_{\eta}}-\ell(\eta)<\lambda}\bigg)\eta^{\ba-\alpha}  \\
& -\widetilde {\mathbb{E}} \bigg( e^{(\ba-\bb)({\widehat
\xi}_{{\widehat T}_{\eta}}-\ell(\eta))}
\Psi(\lambda\!-\!({\widehat \xi}_{{\widehat
T}_{\eta}}\!-\!\ell(\eta)) \nbOne_{{\widehat \xi}_{{\widehat
T}_{\eta}}-\ell(\eta)<\lambda }\bigg).
\end{split}
\end{equation*}
\end{dem}

\medskip

We next summarize the main results of this section in an
asymptotic comparative scheme. The notation $F_{1,2}$ refers to
the situation where in the two-step fragmentation procedure both
devices are effectively used (i.e. $\eta_0/\eta\in (0,1)$),
whereas the notation $F_1$ and $F_2$ respectively refer to the
situations where only the first or only the second fragmentation
process is used.
\begin{cor}\label{summary} Assume that $U(dx)$ and $\widehat{U}(dx)$
have no atom at $0$. In each of the following cases,   the
corresponding assertion holds true for any value of
$\eta_0/\eta\in (0,1)$ as soon as $\eta$ is sufficiently small:
\begin{equation*}
\begin{split}
\widehat{\alpha} >\alpha, \  \widehat{\beta}<\beta
\mbox{ (thus }\alpha - \beta
 < \widehat{\alpha}  -\widehat{\beta})
 :& \ F_1 \mbox{  is better than }F_{1,2}
\mbox{ which is better than }  F_2 \ .  \\
 \widehat{\alpha} <\alpha , \ \widehat{\beta}>\beta
\mbox{ (thus }\alpha - \beta  > \widehat{\alpha}  -\widehat{\beta})
:& \ F_2 \mbox{ is better than } F_{1,2}
\mbox{ which is better than }F_1 \ .   \\
 \widehat{\alpha} <\alpha, \ \widehat{\beta}<\beta
\mbox{ and } \alpha - \beta
 < \widehat{\alpha}  -\widehat{\beta} :& \ F_1
\mbox{ is better than } F_2 \mbox{ which is better than }
\ F_{1,2} \ . \\
\widehat{\alpha} <\alpha, \ \widehat{\beta}<\beta \mbox{ and } \alpha
-\beta  >\widehat{\alpha} -\widehat{\beta}   :& \ F_2
\mbox{ is better than } F_1
\mbox{ which is better than }  \ F_{1,2} \ . \\
\widehat{\alpha} >\alpha, \  \widehat{\beta}>\beta
\mbox{ and }\alpha - \beta <\widehat{\alpha}  -\widehat{\beta} :&
\  F_{1,2} \mbox{ is better than }F_1
\mbox{ which is better than  }F_2 \ . \\
 \widehat{\alpha} >\alpha, \  \widehat{\beta}>\beta  \mbox{ and }
\alpha -
 \beta
 > \widehat{\alpha}  -\widehat{\beta} :& \  F_{1,2}
\mbox{ is better than }F_2 \mbox{ which is better than  }F_1 \ . \\
\end{split}
\end{equation*}
\end{cor}
\smallskip

\begin{rk} By parts c) of  Theorems \ref{2frvs1} and \ref{2frvs2}, if $\widehat{\alpha} =\alpha$  or
if $\widehat{\beta}=\beta$  the comparative efficiency of $F_1$,
$F_2$ and $F_{1,2}$ for $\eta$ small enough is in general
determined by those parameters but  also by the value of
$\eta_0/\eta\in (0,1)$.
\end{rk}

\section{Close-to-unit size thresholds}

We shall next be interested in the behavior of $
\mathbb{E}(\E(\eta ,\eta_0))$ for large values of $\eta$ and
$\eta_0$. Again, we shall compare the mean energy of the two-step
fragmentation procedure with the situations when only the second,
or only the first fragmentation process is used.

\medskip

We shall assume in this analysis  that the subordinators $\xi$ and
${\widehat \xi}$ satisfy under $\widetilde{\mathbb{P}}$
a condition of regular variation at
$\infty$. Namely, respectively denoting by $\widetilde{\phi}$ and
$\widehat{\ttp}$ their Laplace exponents (see
(\ref{phitilde})),  we assume
$$
({\bf RV})\;\;\;\; \exists \; \rho,\, \hpr \, \in \, (0,1) \,
\hbox{ such that } \, \forall \,
\lambda\geq 0 \,: \;\;\; \lim_{q\to \infty}
\frac{\widetilde{\phi}(\lambda q)}{\widetilde{\phi}(q)}=
\lambda^{\rho}\,,\;\;\; \lim_{q\to \infty}
\frac{\widehat{\ttp}(\lambda q)}{\widehat{\ttp}(q)}=
\lambda^{\hpr}\,.
$$

This  assumption can be equivalently (and transparently) stated in
terms of the infinitesimal behavior near $\eta=1$ of the ``mean
energy functions'' $\eta\mapsto \EE({\mE}(\eta))$ and
$\eta\mapsto\EE({\mmE}(\eta))$ of each of the fragmentation
processes. See Remark \ref{avennear0} below.

\medskip

Recall that a function $G:\mathbb{R}_+\to\mathbb{R}_+ $ is said to
vary  slowly  at $0$ if $\lim_{x\to 0^+} G(\lambda x) /G(x) =1$
for all $\lambda \in [0,\infty)$. A well known fact that will be
used in the sequel is that such convergence is uniform in $\lambda
\in [0,\lambda_0]$, for all $\lambda_0 \in (0,\infty)$.

\medskip

By $L$ and $\hL$ we shall denote the nonnegative slowly varying
functions at $0$ defined by the relations
$$
L\left(\frac{1}{ x}\right) = \frac{x^{\rho}}{\widetilde{\phi}(x)}
\,,\;\;\; {\hL}\left(\frac{1}{x}\right)=
\frac{x^{\hpr}}{\widehat{\ttp}(x)}.
$$

\begin{rk}
Using the aforementioned uniform convergence result for $L$ and
$\hL$ it is not hard to check  that $ \lim\limits_{q \to
\infty}\frac{{\phi}( q)}{\widetilde{\phi}(q)}= \lim\limits_{q \to
\infty}\frac{{\hp}( q)}{\widehat{\ttp}(q)} =1$. Consequently,
({\bf RV}) implies that the same condition hold on $\phi$ and
$\hp$ and conversely.
\end{rk}

We  define
$$
Q_{\phi,\hp}:=\lim\limits_{q \to \infty}
\frac{\hp( q)}{\phi(q)} \ = \lim\limits_{q \to
\infty}\frac{\widehat{\ttp}(q)}{\widetilde{\phi}(q)} =
\lim\limits_{x \to 0^+} \frac{L(x)x^{\rho}}{\hL(x)x^{\hpr}}
$$
if the limit in $ [0,\infty]$ exists.
More generally, we write
$$
Q_{{\phi},\hp}^+ :=\limsup\limits_{q \to \infty}
\frac{\hp( q)}{\phi(q)} \ = \limsup\limits_{x \to 0^+}
\frac{L(x)x^{\rho}}{\hL(x)x^{\hpr}}
$$
and
$$
Q_{{\phi},{\hp}}^- :=\liminf\limits_{q \to
\infty}\frac{\hp( q)}{\phi(q)} \ = \liminf\limits_{x \to 0^+}
\frac{L(x)x^{\rho}}{\hL(x)x^{\hpr}}.
$$

Recall the notation
$$
\Psi(x)=C \, \int_{0}^{x} e^{(\alpha-\beta)y}
{U} (dy)\,,\;\;\; \hPs(x)=\cc \,\int_{0}^{x}
e^{(\ba-\bb) y} \widehat{U} (dy) .
$$

\begin{lm}
\label{firstpartofthelimit}
We have
$$
C Q_{{\phi},{\hp}}^- -\cc\leq  \liminf\limits_{\eta\to 1^-}
\frac{\Psi(\ell(\eta))
-\hPs(\ell(\eta))}{\hL(\ell(\eta)) {\ell(\eta)}^{\hpr}} \leq
\limsup\limits_{\eta\to 1^-} \frac{\Psi(\ell(\eta))
-\hPs(\ell(\eta))}{\hL(\ell(\eta)) {\ell(\eta)}^{\hpr}}\leq  C
Q_{{\phi},{\hp}}^+ -\cc .
$$
In particular,
\begin{equation*}
\lim\limits_{\eta\to 1^-} \frac{\Psi(\ell(\eta))
-\hPs(\ell(\eta))}{\hL(\ell(\eta))
{\ell(\eta)}^{\hpr}}= \begin{cases} \infty
& \mbox{ if } \hpr>\rho \\
 - \cc & \mbox{ if } \hpr<\rho\\ C Q_{{\phi},{\hp}} -\cc & \mbox{ if }
\hpr=\rho \mbox{ and } \exists \;
Q_{{\phi},{\hp}}=\lim\limits_{x \to 0^+} \frac{L(x)}{\hL(x)}
\in [0,\infty]\,.
\end{cases}
\end{equation*}

\end{lm}

\begin{dem} By classic Tauberian theorems (see e.g. Th, 5.13 in \cite{Ky}
or Section 0.7 in \cite{be1}), our assumptions on
$\widetilde{\phi}$ and $\widehat{\ttp}$ are respectively
equivalent to
\begin{equation*}
\lim\limits_{x \to 0^+} \frac{{U}(x)}{x^{\rho}L(x)}=1
\,,\;\;\; \lim\limits_{x \to 0^+}
\frac{\widehat{U}(x)}{x^{\hpr}\hL(x)}= 1\,.
\end{equation*}
On the other hand, we have
\begin{equation*}
\begin{split}
\Psi(x) -\hPs(x)\leq & \ C e^{|\alpha-\beta|x}
{U}(x)-\cc e^{-|\ba-\bb|x}\widehat{U}(x) \\
 = & \hL(x)x^{\hpr} \cc \left( \frac{C}{\cc} e^{|\alpha-\beta|x}
\frac{{U}(x)}{x^{\rho}L(x)}
\frac{L(x)}{\hL(x)}x^{\rho-\hpr} -e^{-|\ba-\bb|x}
\frac{\widehat{U}(x)}{x^{\hpr}\hL(x)} \right)\\
\end{split}
\end{equation*}
and similarly,
\begin{equation*}
\begin{split}
\Psi(x)-\hPs(x)\geq & \ C e^{-|\alpha-\beta|x}
{U}(x)-\cc e^{|\ba-\bb|x}\widehat{U}(x) \\
 = & \hL(x)x^{\hpr} \cc \left( \frac{C}{\cc} e^{-|\alpha-\beta|x}
\frac{{U}(x)}{x^{\rho}L(x)}
\frac{L(x)}{\hL(x)}x^{\rho-\hpr} -e^{|\ba-\bb|x}
\frac{\widehat{U}(x)}{x^{\hpr}\hL(x)} \right)\,.\\
\end{split}
\end{equation*}
The first statement follows from these bounds. To complete the
proof, notice that since $\frac{L(x)}{\hL(x)}$  is slowly varying
at $0$, we have that
\begin{equation*}
\lim_{x\to 0^+}
\frac{L(x)}{\hL(x)}x^{\rho-\hpr}=
\begin{cases} \infty
& \mbox{ if }\hpr>\rho \\
 0 & \mbox{ if } \hpr<\rho\\  Q_{{\phi},{\hp}}  & \mbox{ if }
\hpr=\rho \mbox{ and } \exists \; \lim\limits_{x \to 0^+}
\frac{L(x)}{\hL(x)}\in [0,\infty]\,,
\end{cases}
\end{equation*}
using also the fact that $\lim\limits_{x\to 0^+}G(x)=0$ for any
 regularly varying (at $0$) function $G(x)$ with positive index.
\end{dem}

\medskip

Notice that ({\bf RV})  implies that $U$ has no atom at $0$ (see
e.g. the first lines of the previous proof).

\begin{rk}
\label{avennear0}
The estimates used in the proof of Lemma
\ref{firstpartofthelimit} show that $$\Psi(x)\sim C
{U}(x)\quad \mbox{ and }\quad \hPs(x)\sim \cc
\widehat{U}(x)\quad \mbox{ when }x\to 0^+,$$ so that
$\Psi(x)\sim C x^{\rho}L(x)$ and $\hPs(x)\sim \cc
x^{\hpr}\hL(x)$ as well. Consequently, by the aforementioned
Tauberian results, assumption ({\bf RV}) is equivalent to
\medskip

({\bf RV}) $x\mapsto\EE({\mE}(e^{-x}))$  and
$x\mapsto\EE({\mmE}(e^{-x}))$ are regularly varying at $0^+$
with indexes $\rho,\, \hpr\in (0,1)$ respectively.

\medskip

This alternative formulation  has the advantage of providing a way
to infer the regularity indexes from separate observations of both
fragmentation processes, if one was able  to measure the energies
required to obtain fragments of different  close to unit sizes.
More precisely,
$$
\frac{\log  \EE({\mE}(\eta^{\lambda}))-
\log \EE({\mE}(\eta))}{\log \lambda}
$$
should be close to $\rho$ for $\eta$ sufficiently close to $1$.
Alternatively, $\rho$  could in principle also be deduced from the
estimation method of $\phi$ developed in \cite{kre}.
\medskip

In the same vein, we remark that   the existence of the limit
$Q_{\phi,\hp}$ is equivalent to
$$
\exists \; Q:=\lim\limits_{\eta\to 1^-}
\frac{\EE({\mE}(\eta))}{\EE({\mmE}(\eta))}=\lim\limits_{\eta\to
1^-}\frac{C}{\cc}Q_{\phi,\hp}.
$$
In general, Lemma \ref{firstpartofthelimit} indeed shows  that
$$
\cc(Q^- -1)\leq  \liminf\limits_{\eta\to 1^-} \frac{\Psi(\ell(\eta))
-\hPs(\ell(\eta))}{\hL(\ell(\eta)) {\ell(\eta)}^{\hpr}} \leq
\limsup\limits_{\eta\to 1^-} \frac{\Psi(\ell(\eta))
-\hPs(\ell(\eta))}{\hL(\ell(\eta)) {\ell(\eta)}^{\hpr}}\leq
\cc(Q^+ -1),
$$
where
$$
Q^+:=\limsup\limits_{\eta\to 1^-}
\frac{\EE({\mE}(\eta))}{\EE({\mmE}(\eta))}=\frac{C}{\cc}Q_{\phi,\hp}^+ \;,
$$
and
$$
Q^-:=\liminf\limits_{\eta\to 1^-}
\frac{\EE({\mE}(\eta))}{\EE({\mmE}(\eta))}=\frac{C}{\cc}
Q_{\phi,\hp}^-\;.
$$

\end{rk}

We recall now that, under our assumptions on the Laplace exponents
$\widetilde{\phi}$ and $\widehat{\ttp}$, by the
Dynkin-Lamperti Theorem it weakly holds as $\eta\to 1^-$ that
$$
\widetilde{\mathbb{P}}
\left(\frac{\xi_{T_{\eta}}-\ell(\eta)}{\ell(\eta)}\in dy\right)
\to \mu(dy):=\frac{\sin(\rho \pi)}{\pi} \frac{dy}{(1+y)y^{\rho}}
$$
and
$$
\widetilde{\mathbb{P}} \left(\frac{{\widehat \xi}_{{\widehat
T}_{\eta}}-\ell(\eta)}{\ell(\eta)}\in dy\right)\to \hmu(dy):=
\frac{\sin(\hpr \pi)}{\pi} \frac{dy}{(1+y)y^{\hpr}}\,.
$$

This  suggest us the way in which $\eta$ and $\eta_0$ should go to
$1$ in order to observe a coherent  close-to-unit size asymptotic
behavior. In all the sequel $\gamma>1$ is a fixed parameter, and
we assume that
$$
\eta_0=\eta_0(\eta)=\eta^{\gamma}.
$$

We have the following

\begin{lm}\label{computationofseconterm}
\begin{multline}
\label{secondterm}
\lim_{\eta\to 1^-}\frac{\widetilde {\mathbb{E}}
\left( \nbOne_{\xi_{T_{\eta}}< \ell(\eta_0)}e^{(\alpha-\bb)
\xi_{T_{\eta}}}\hPs(\ell(\eta_0)-\xi_{T_{\eta}})\right)
 -\widetilde{\mathbb{E}}
\left(\nbOne_{{\widehat \xi}_{{\widehat T}_{\eta}}<
 \ell(\eta_0)}e^{(\ba-\bb){\widehat \xi}_{{\widehat T}_{\eta}}}
\hPs(\ell(\eta_0)-{\widehat \xi}_{{\widehat T}_{\eta}})\right)}
{(\ell(\eta))^{\hpr}\hL(\ell(\eta))} \\
= \cc \left[ \int_0^{\gamma-1}\!\!\! (\gamma-1-y)^{\hpr}\mu(dy)-
\int_0^{\gamma-1} \!\!\! (\gamma-1-y)^{\hpr} \hmu(dy)\right]\,.
\end{multline}
Moreover, in the case $\frac{1}{2}\geq \rho>\hpr$,  the limit is
a nonnegative and increasing  function of $\gamma$ for $\gamma\in
[1,2]$, which goes to $0$ when $\gamma \to 1^+$.

\end{lm}

\begin{dem} Denote by $\partial(\eta)$ the numerator in the left hand
side \eqref{secondterm} and respectively by $\mu^{\eta}$ and
${\hmu}^{\eta}$ the laws of
$$
\frac{\xi_{T_{\eta}}-\ell(\eta)}{\ell(\eta)} \;\, \hbox{ and } \;\,
\frac{{\widehat \xi}_{T_{\eta}}-\ell(\eta)}{\ell(\eta)}.
$$
We then easily see that
\begin{equation*}
\begin{split}
\partial(\eta) \leq & e^{\gamma \ell(\eta)
|\alpha-\bb|}
\int_0^{\gamma-1}\hps(\ell(\eta)(\gamma-1-y))\mu^{\eta}(dy)
\\& - e^{- \gamma \ell(\eta) |\ba-\bb|}
\int_0^{\gamma-1}\hps(\ell(\eta)(\gamma-1-y)){\hmu}^{\eta}(dy)\\
\end{split}
\end{equation*}
and
\begin{equation*}
\begin{split}
\partial(\eta) \geq & e^{-\gamma \ell(\eta)  |\alpha-\bb|}
\int_0^{\gamma-1}\hps(\ell(\eta)(\gamma-1-y))\mu^{\eta}(dy)
\\& - e^{\gamma \ell(\eta) |\ba-\bb|}
\int_0^{\gamma-1}\hps(\ell(\eta)(\gamma-1-y)){\hmu}^{\eta}(dy).\\
\end{split}
\end{equation*}
On the other hand, by similar estimates as in the previous lemma, one
checks that
\begin{equation}
\label{theta} \theta(x):=\frac{\hps(x)}{\cc \hL(x)x^{\hpr}}\to 1
\end{equation}
when $x\searrow 0$, and thus $\theta(x)$ is slowly varying at
$0$. Fix now $\varepsilon\in (0,1)$, and recall that for a slowly varying
at $0$ function $G(x)$, the convergence $G(\lambda x)/G(x)\to 1$
is uniform in $\lambda \in [0,\lambda_0]$ for all $\lambda_0\in (0,1)$.
Therefore, since
\begin{equation*}
\hps(\ell(\eta) y)=
\frac{\theta(\ell(\eta)y)}{\theta(\ell(\eta))}
\frac{\hL(\ell(\eta) y)}{\hL(\ell(\eta))}
\hps(\ell(\eta))y^{\hpr},
\end{equation*}
we deduce that if $\eta\in (0,1)$ is sufficiently close to $1$,
\begin{equation*}
\forall \, y\in [0,\gamma-1]:\;\;\;
(1-\varepsilon)\hps(\ell(\eta))y^{\hpr} \leq \hps(\ell(\eta) y )
\leq (1+\varepsilon)\hps(\ell(\eta)) y ^{\hpr}\,.
\end{equation*}
Moreover, from \eqref{theta}, it follows that
if $\eta$ is sufficiently close to $1$ then
\begin{equation}
\label{ineqs}
\forall \, y\in [0,\gamma-1] :\;\;\;
\cc (1-\varepsilon)^2 y^{\hpr} \leq \frac{\hps(\ell(\eta) y )}
{(\ell(\eta))^{\hpr}\hL(\ell(\eta))}
\\ \leq \cc (1+\varepsilon)^2 y ^{\hpr}\,.
\end{equation}
It follows that
\begin{equation*}
\limsup_{\eta\to 1^-} \frac{\partial(\eta)}{{(\ell(\eta))^{\hpr}
\hL(\ell(\eta))}} \leq (1+\varepsilon)^2 \cc A_{\gamma}-
(1-\varepsilon)^2 \cc {\ca}_{\gamma}
\end{equation*}
and
\begin{equation*}
\liminf_{\eta\to 1^-}
\frac{\partial(\eta)}{{(\ell(\eta))^{\hpr} \hL(\ell(\eta))}}
\geq (1-\varepsilon)^2 \cc A_{\gamma}- (1+\varepsilon)^2
\cc {\ca}_{\gamma},
\end{equation*}
where
$$
A_{\gamma}=\frac{\sin(\pi\rho)}{\pi}\!\int_0^{\gamma-1} \!\!\!\!
(\gamma\!-\!1\!-\!u)^{\hpr} \frac{du}{(1+u)u^{\rho}}\,, \;\;
{\ca}_{\gamma}=\frac{ \sin(\pi \hpr)}{\pi} \! \int_0^{\gamma-1}
\!\!\!\! (\gamma\!-\!1\!-\!u)^{\hpr} \frac{du}{(1+u)u^{\hpr}} \,.
$$
The first statement  follows by letting $\varepsilon\to 0^+$.  The
asserted properties of $\cc(A_{\gamma}-{\widehat A}_{\gamma})$ are
consequence of the inequalities $u^{-\rho}>u^{-\hpr}$ for $u\in
(0,1)$, $\sin(\pi \rho)> \sin(\pi \hpr)>0$ when $\frac{1}{2}>\rho>
\hpr$, and dominated convergence.
\end{dem}

\medskip

We next introduce  helpful concepts in order to state our results
on the energy for large thresholds.

\begin{df}
\noindent $(i)$
The fragmentation processes $\X$ is said to be infinitesimally
efficient ({\it inf. eff.}) compared to $\hX$ if ({\bf RV}) holds  and
$Q_{\phi,\hp}^+<\frac{\cc}{C}.$

\medskip

Conversely,

\medskip

\noindent $(ii)$ The fragmentation processes $\hX$ is said to be
{\it inf. eff.} compared to $\X$  if ({\bf RV}) holds
and $Q_{\phi,\hp}^->\frac{\cc}{C}.$
\end{df}

For instance,  $\X$ is {\it inf. eff.} compared to
$\hX$ if  $\rho> \hpr$ or if $\rho=\hpr$ and $Q_{\phi,\hp}$
exists in $[0,\frac{\cc}{C})$. Similarly, $\hX$ is
{\it inf. eff.} compared to $\hX$ e.g. if $\rho< \hpr$
or if $\rho=\hpr$ and $Q_{\phi,\hp}$ exists in
$(\frac{\cc}{C},\infty]$.

\medskip

\begin{rk}
We observe that  $\X$ (respectively $\hX$) is {\it inf. eff.}
compared to $\hX$ (respectively $\X$) if and only if
$\EE({\mE}(e^{-x}))$ and $\EE({\mmE}(e^{-x}))$ are
regularly varying functions at $x=0$ with indexes in $(0,1)$ and
$Q^+=\limsup\limits_{\eta\to 1^-}
\frac{\EE({\mE}(\eta))}{\EE({\mmE}(\eta))}<1$ (respectively
$Q^-=\liminf\limits_{\eta\to 1^-}
\frac{\EE({\mE}(\eta))}{\EE({\mmE}(\eta))}>1$).
\end{rk}
\medskip

Bringing all together, we have obtain:

\begin{theo}
\label{mainforetacloseto1} (Two-step procedure versus second
fragmentation only)

For each $\gamma \in (1,\infty)$ it holds:

\smallskip

\noindent $(a)$ If $\hX$ is {\it inf. eff.} compared to
$\X$ and $Q^-=Q=\infty$ (in particular if $\hpr >\rho$), then:

$\forall \, M> 0\;$ $\exists\, \eta_M\in (0,1)$ such that
\begin{equation*}
\forall \, \eta\in (\eta_M,1]: \;\;\;
\mathbb{E}(\E(\eta ,\eta^{\gamma}))>
\mathbb{E}({\mmE}(\eta^{\gamma}))+M
\mathbb{E}({\mmE}(\eta))>
\mathbb{E}({\mmE}(\eta^{\gamma})).
\end{equation*}

\noindent $(b)$ If $\hX$ is {\it inf. eff.} compared to
$\X$ and $Q^-\in (1,\infty)$ (and thus $\rho=\hpr$), then:

$\forall \, \varepsilon\in (0,Q^- -1)\;$
$\exists\, \eta_{\varepsilon}\in (0,1)$ such that
\begin{equation*}
\forall \, \eta\in  (\eta_{\varepsilon},1] : \;\;\;
\mathbb{E}(\E(\eta ,\eta^{\gamma}))>
\mathbb{E}({\mmE}(\eta^{\gamma}))+ (Q^-
-1-\varepsilon)\mathbb{E}({\mmE}(\eta)) >
\mathbb{E}({\mmE}(\eta^{\gamma})).
\end{equation*}

\noindent $(c)$ If $\X$ is {\it inf. eff.} compared to
$\hX$ and $Q^+\in (0,1)$  (and thus $\rho=\hpr$), then:

$\forall \, \varepsilon\in (0,1-Q^+)\;$
$\exists\, \eta_{\varepsilon}\in (0,1)$ such that
\begin{equation*}
\forall \, \eta\in (\eta_{\varepsilon},1]: \;\;\;
\mathbb{E}(\E(\eta ,\eta^{\gamma})) <
\mathbb{E}({\mmE}(\eta^{\gamma}))+
(Q^+-1+\varepsilon)\mathbb{E}({\mmE}(\eta))<\mathbb{E}({\mmE}(\eta^{\gamma})).
\end{equation*}

\noindent $(d)$ If $\X$ is {\it inf. eff.} compared to
$\hX$ and $Q^+=Q=0$ (in particular if $\hpr<\rho$), then:

$\forall \, \varepsilon\in (0,1)$, $\exists\, \eta_{\varepsilon}\in (0,1)$
such $\forall \, \eta\in (\eta_{\varepsilon},1]$:
\begin{multline*}
\mathbb{E}({\mmE}(\eta^{\gamma}))+
(A_{\gamma}-{\ca}_{\gamma}-1-\varepsilon)\mathbb{E}({\mmE}(\eta))
\!< \!\mathbb{E}(\E(\eta ,\eta^{\gamma}))\!< \!
\mathbb{E}({\mmE}(\eta^{\gamma}))\!+\!
(A_{\gamma}-{\ca}_{\gamma}-1+\varepsilon)\mathbb{E}({\mmE}(\eta)).
\end{multline*}
(The quantities $A_{\gamma}$ and ${\ca}_{\gamma}$ were defined in
Lemma \ref{computationofseconterm}).

\medskip

Moreover, if $\frac{1}{2}\geq   \rho>\hpr$,
$\exists \, \gamma_0\in (1,2]$ such that
$\forall \, \gamma \in (1,\gamma_0]$, one has
$1-A_{\gamma}+{\ca}_{\gamma}>0$ and
$\forall \, \varepsilon \in (0,1-A_{\gamma}+{\ca}_{\gamma})$,
\begin{equation*}
\forall \, \eta\in (\eta_{\varepsilon},1]: \;\;\;
\mathbb{E}(\E(\eta ,\eta^{\gamma}))<
\mathbb{E}({\mmE}(\eta^{\gamma}))+
(A_{\gamma}-{\ca}_{\gamma}-1+\varepsilon)
\mathbb{E}({\mmE}(\eta))<\mathbb{E}({\mmE}(\eta^{\gamma})).
\end{equation*}

In all four cases,  similar statements hold with
$\mathbb{E}({\mmE}(\eta))$ replaced by
$\cc
\left[\widehat{\ttp}\left(\frac{1}{\log (1/\eta)}\right)\right]^{-1}$.
\end{theo}
\begin{dem}
By Remark \ref{decomp} and the previous results, we simply have to
notice that when $\eta \to 1^-$,
$$
\frac{\EE({\mmE}(\eta))}{\cc}\sim
\hL(\ell(\eta))(\ell(\eta))^{\hpr} =
\left[\widehat{\ttp}\left(\frac{1}{\ell(\eta)}\right)\right]^{-1},
$$
$\E(1 ,\eta^{\gamma})= {\mmE}(\eta^{\gamma})$ and
the quantities $A_{\gamma}$ and ${\ca}_{\gamma}$ are equal if $\rho=\hpr$.
The last assertion in part $(d)$ is  consequence of the last
part of Lemma \ref{computationofseconterm}.
\end{dem}

\medskip

The previous theorem provided conditions on large thresholds
$\eta$ and $\eta_0$ under which the use of the second
fragmentation process can be told to be efficient or not. We next
briefly address  the efficiency of using or not the first
fragmentation process. The arguments of the following theorem are
similar to those of the previous lemmas, so we  just sketch its
proof. We use the following notation
$$
\forall \, \gamma \in (1,\infty): \;\;\;
B_{\gamma}:=\frac{\sin(\pi\rho)}{\pi} \! \int_0^{\gamma-1}
\!\!(\gamma\!-\!1\!-\!u)^{\rho} \frac{du}{(1+u)u^{\rho}}\,.
$$

\begin{theo} (Two-step procedure versus first
fragmentation only)

For all $\gamma \in (1,\infty)$ it holds:

\medskip

\noindent $(a)$ If $\hX$ is {\it inf. eff.} compared to
$\X$, then:

$\forall \, \varepsilon\in (0,1-\frac{1}{Q^-})\;$
$\exists \, \eta_{\varepsilon}\in (0,1)$ such that
\begin{equation*}
\forall \, \eta\in (\eta_{\varepsilon},1]: \;\;\;
\mathbb{E}(\E(\eta ,\eta^{\gamma}))<
\mathbb{E}({\mE}(\eta^{\gamma}))+
\left(\frac{1}{Q^-}-1 +\varepsilon\right)
B_{\gamma}\mathbb{E}({\mE}(\eta)) <
\mathbb{E}({\mE}(\eta^{\gamma})).
\end{equation*}

\noindent $(b)$  If $\X$ is {\it inf. eff.} compared to
$\hX$ and $Q^+\in (0,1)$  (and thus $\rho=\hpr$), then:

$\forall \, \varepsilon\in (0,\frac{1}{Q^+}-1)\;$,
$\exists \, \eta_{\varepsilon}\in (0,1)$ such that
\begin{equation*}
\forall \, \eta\in (\eta_{\varepsilon},1]: \;\;\;
\mathbb{E}(\E(\eta ,\eta^{\gamma}))>
\mathbb{E}({\mE}(\eta^{\gamma}))+
\left(\frac{1}{Q^+}-1 -\varepsilon\right)
B_{\gamma}\mathbb{E}({\mE}(\eta)) >
\mathbb{E}({\mE}(\eta^{\gamma})).
\end{equation*}

\noindent $(c)$ If $\X$ is {\it inf. eff.} compared to
$\hX$ and $Q^+=Q=0$ (in particular if $\hpr<\rho$), then:

$\forall \, M> 0\;$ $\exists \, \eta_M\in (0,1)$ such that
\begin{equation*}
\forall \, \eta\in (\eta_M,1]: \;\;\;
\mathbb{E}(\E(\eta ,\eta^{\gamma}))>
\mathbb{E}({\mE}(\eta^{\gamma}))+M \mathbb{E}({\mE}(\eta))>
\mathbb{E}({\mE}(\eta^{\gamma})).
\end{equation*}

In all cases, one can replace $\mathbb{E}({\mE}(\eta))$  by
$C\left[\widetilde{\phi}\left(\frac{1}{\ell(\eta)}\right)\right]^{-1}$.

\end{theo}
\begin{dem}
Fix $\gamma >1$ and $\varepsilon\in (0,1)$.  As in Lemma
\ref{computationofseconterm} we get that for all $y\in
[0,\gamma-1]$,
\begin{equation*}
C (1-\varepsilon)^2 y^{\rho} \leq \frac{\psi(\ell(\eta) y )}
{(\ell(\eta))^{\rho}L(\ell(\eta))}
\\ \leq C (1+\varepsilon)^2 y ^{\rho}
\end{equation*}
and
\begin{equation*}
\cc (1-\varepsilon)^2
\frac{(\ell(\eta))^{\hpr}\hL(\ell(\eta))}{(\ell(\eta))^{\rho}L(\ell(\eta))}
y^{\hpr} \leq \frac{\hps(\ell(\eta) y )}
{(\ell(\eta))^{\rho}L(\ell(\eta))}\\
\leq \cc (1+\varepsilon)^2  \frac{(\ell(\eta))^{\hpr}
\hL(\ell(\eta))}{(\ell(\eta))^{\rho}L(\ell(\eta))}
y ^{\hpr}
\end{equation*}
if  $\eta$ is close enough to $1$. Set now
$\overline{\partial}(\eta):=\mathbb{E}(\E(\eta
,\eta^{\gamma})-\E(\eta^{\gamma},\eta^{\gamma}))$. From the
previous bounds,  and from  the explicit expression for
$\overline{\partial}(\eta)$ given in Remark \ref{decomp}, we
deduce  that
\begin{equation*}
C \left(\frac{A_{\gamma}}{Q^+}- B_{\gamma}\right)
\leq \liminf_{\eta\to 1^-}
\frac{\overline{{\partial}}(\eta)}{{(\ell(\eta))^{\rho}
L(\ell(\eta))}} \leq \limsup_{\eta\to 1^-}
\frac{\overline{\partial}(\eta)}{{(\ell(\eta))^{\rho}
L(\ell(\eta))}} \leq C\left(\frac{A_{\gamma}}{Q^-} - B_{\gamma}\right).
\end{equation*}
Part  $(a)$ follows from this relation, using the facts that
$A_{\gamma}=B_{\gamma}$ if $\rho =\hpr$, and that $Q_- =\infty$
if $\hpr> \rho$. The remaining parts are similar.
\end{dem}

\medskip

\begin{rk}
If $\rho=\hpr$ and $Q\in (0,\infty)$ exists,  one obtains
for $\eta$ close enough to $1$,
\begin{equation*}
\mathbb{E}({\mmE}(\eta^{\gamma}))+
(Q-1-\varepsilon)\mathbb{E}({\mmE}(\eta))  <
\mathbb{E}(\E(\eta ,\eta^{\gamma}))<
\mathbb{E}({\mmE}(\eta^{\gamma}))+
(Q-1+\varepsilon)\mathbb{E}({\mmE}(\eta))
\end{equation*}
\begin{equation*}
\mathbb{E}({\mE}(\eta^{\gamma}))+
(Q^{-1}-1-\varepsilon)B_{\gamma}\mathbb{E}({\mE}(\eta))  <
\mathbb{E}(\E(\eta ,\eta^{\gamma}))<
\mathbb{E}({\mE}(\eta^{\gamma}))+
(Q^{-1}-1+\varepsilon)B_{\gamma}\mathbb{E}({\mE}(\eta)).
\end{equation*}
In particular, when $Q=1$ we deduce that
$\mathbb{E}(\E(\eta ,\eta^{\gamma}))\sim
\mathbb{E}({\mmE}(\eta^{\gamma}))\sim
\mathbb{E}({\mE}(\eta^{\gamma}))$ when $\eta\to 1^-$, as one
could expect.
\end{rk}

\noindent {\bf Acknowledgments.} J.Fontbona and S.Mart\'\i nez are
indebted to Basal Conicyt Project.

\medskip

\noindent JOAQU\'IN FONTBONA

\noindent {\it Departamento Ingenier{\'\i}a Matem\'atica and Centro
Modelamiento Matem\'atico, Universidad de Chile,
UMI 2807 CNRS, Casilla 170-3, Correo 3, Santiago, Chile.}

e-mail: $<fontbona@dim.uchile.cl>$

\medskip

\noindent NATHALIE KRELL

\noindent  {\it IRMAR,
Universit\'e Rennes 1,
Campus de Beaulieu,
35042 Rennes Cedex FRANCE.}

e-mail: $<nathalie.krell@univ-rennes1.fr>$

\medskip

\noindent SERVET MART\'INEZ

\noindent {\it Departamento Ingenier{\'\i}a Matem\'atica and Centro
Modelamiento Matem\'atico, Universidad de Chile, UMI 2807 CNRS,
Casilla 170-3, Correo 3, Santiago, Chile.}

e-mail: $<smartine@dim.uchile.cl>$

\end{document}